\documentclass[11pt,oneside]{amsart}
\usepackage{amsthm,enumitem}
\usepackage[
    a4paper,
    left=2.5cm,
    right=2.5cm,
    top=2.5cm,
    bottom=2.5cm,
    marginparwidth=3cm,
    marginparsep=0.5cm
]{geometry}

\usepackage{amsmath,amsfonts,amssymb,graphicx} 
\usepackage{cancel}

\usepackage{mathpazo} 
\usepackage[dvipsnames]{xcolor} 
\usepackage{hyperref}
\usepackage[normalem]{ulem}

\usepackage{comment}

\usepackage{amsthm}
\usepackage{aliascnt}
\usepackage{hyperref}

\newtheorem{thm}{Theorem}

\newaliascnt{cor}{thm}
\newtheorem{cor}[cor]{Corollary}
\aliascntresetthe{cor}

\newaliascnt{lem}{thm}
\newtheorem{lem}[lem]{Lemma}
\aliascntresetthe{lem}

\newaliascnt{prop}{thm}
\newtheorem{prop}[prop]{Proposition}
\aliascntresetthe{prop}

\newaliascnt{defn}{thm}

\aliascntresetthe{defn}

\newaliascnt{ex}{thm}

\aliascntresetthe{ex}

\newaliascnt{rem}{thm}
\newtheorem{rem}[rem]{Remark}
\aliascntresetthe{rem}

\newaliascnt{conj}{thm}

\aliascntresetthe{conj}

\newaliascnt{qu}{thm}

\aliascntresetthe{qu}

\newcommand{\thistheoremname}{}
\newtheorem*{genericthm*}{\thistheoremname}

\newcounter{namedthm}

\newenvironment{namedthm*}[1]
  {%
    \renewcommand{\thistheoremname}{#1}%
    \refstepcounter{namedthm}%
    \begin{genericthm*}%
  }
  {%
    \end{genericthm*}%
  }

\newcommand{\N}{\mathbb{N}}

\newcommand{\R}{\mathbb{R}}
\newcommand{\C}{\mathbb{C}}

\newcommand{\Hc}{\mathcal{H}}

\newcommand{\Ac}{\mathcal{A}}

\newcommand{\rmC}{\mathrm{C}}

\DeclareMathOperator{\tr}{tr}

\def\multinom#1#2{\ensuremath{\left(\kern-.3em\left(\genfrac{}{}{0pt}{}{#1}{#2}\right)\kern-.3em\right)}}

\author[]{Oliver F\"urst}
\address{Mathematisches Institut, Universit\"at Bonn, Endenicher Allee 60, 53115 Bonn, Germany.}
\email{ofuerst@math.uni-bonn.de}

\author[]{Eva-Maria Hekkelman}
\address{Max Planck Institute for Mathematics, Vivatsgasse 7, 53111 Bonn, Germany.}
\email{hekkelman@mpim-bonn.mpg.de}

\author[]{Malte Leimbach}
\address{Christian-Albrechts-Universit\"at zu Kiel, Mathematisches Seminar, Heinrich-Hecht-Platz 6, 24118 Kiel, Germany.}
\email{leimbach@math.uni-kiel.de}

\title{Generic properties of the Laplace spectrum and eigenfunctions}

\begin{document}

\date{29 October, 2026}

\maketitle

\begin{abstract}
    We prove that the spectrum of the Laplace--Beltrami operator on compact Riemannian manifolds without boundary is generically linearly independent over $\mathbb{Q}$. Furthermore, we prove that pointwise products of eigenfunctions generically are linearly independent. The proof is based on conformal perturbations of the metric. The result implies in particular that generically the truncation of the continuous functions on $M$ by spectral projections of the Laplace operator is given by symmetric matrices.
\end{abstract}

\section{Introduction}
In contrast with the extremely extensive literature on the spectrum of Laplace operators on manifolds, examples for which we can compute this spectrum explicitly are exceedingly rare. 
Notably, on the $50$th anniversary of Uhlenbeck's celebrated result showing that the spectrum of the Laplace--Beltrami operator on a generic manifold is simple~\cite{Uhlenbeck1976}  (meaning that all eigenspaces are $1$-dimensional),  not a single explicit example exhibiting this behaviour has been produced: in a hay stack, we have only found needles. In fact, this illustrates the need for techniques studying generic properties of spectra and eigenfunctions, as the computable examples can be highly misleading compared to the bigger picture.

Uhlenbeck's genericity result for spectra has since been reproduced for many different classes of examples. 
Where she proved  generic simplicity of the Laplace--Beltrami spectrum  for $C^k$-metrics,  the same property  has now been proven for the Laplace--Beltrami spectrum for smooth Riemannian metrics~\cite{Tanikawa1979, BleeckerWilson1980, BandoUrakawa1983}, the spin Dirac spectrum (up to algebraically forced multiplicities) on manifolds of dimension $2$ and $3$~\cite{Dahl2003}, and most recently, for perturbations of sub-Riemannian Laplace operators and maximally hypoelliptic pseudodifferential operators by bounded potentials~\cite{ChatzakouHelffer2026}.

In this paper, we extend Uhlenbeck's result in a different direction. 
We prove that for a generic  connected closed (i.e.\ compact without boundary)  smooth Riemannian $(M,g)$  of dimension at least $2$,  the spectrum of a Laplace--Beltrami operator $\sigma(\Delta_g) =\{\lambda_i\}_{i \in \mathbb{N}}$ (always listed in increasing order) exhibits no extra hidden structure; it is linearly independent over the rational numbers. 
Here we really mean the spectrum without the first eigenvalue $\lambda_0=0$, that is, for all $q_1,\ldots,q_N\in \mathbb{Q}$, we prove that
\[
\sum_{i=1}^N q_{i} \lambda_i = 0\quad  \Rightarrow \quad  q_1=\cdots=q_N=0.
\]
Observe that this implies that the spectrum is simple. 
Loosely speaking,
in the tradition of \cite{kac66-drum}, this furthermore shows that the sound a drum generically contains no harmonies whatsoever. 
For this reason, one might subtitle our result as 
\begin{quote}
\begin{center}
    \emph{Almost all drums sound awful}.
\end{center}
\end{quote}
Moreover, for the corresponding eigenfunctions $\{e_i\}_{i \in \mathbb{N}}$ of the Laplace--Beltrami operator 
we prove that generically, their pointwise products are linearly independent (over $\mathbb{C}$). For a finitely supported multi-index $\alpha=(\alpha_i)_{i\ge1}\in\mathbb N^{(\infty)}$,
we write
\[
  e^\alpha := \prod_{i\ge 1} e_i^{\alpha_i},\qquad
  |\alpha| := \sum_{i\ge1}\alpha_i,
\]
with the convention $e^0 = e_0$. 

Now to make the word `generic' more precise, recall that a residual set in a topological space is a set which can be written as a countable intersection of open dense sets. 
 For a connected closed smooth manifold, we denote by $\mathcal{R}(M)$ the set of Riemannian metrics on $M$. 
This has the structure of a Fr\'echet manifold: 
a sequence of metrics $\{g^{n}\}_{n\in \N}$ converges if in each chart of $M$ the smooth functions $\{g_{ij}^{n}(x)\}_{n \in \N}$ converge in the usual Fr\'echet topology for smooth functions. 
See \cite{BandoUrakawa1983} for the metric space structure on $\mathcal{R}(M)$.

For any given metric $g \in \mathcal{R}(M)$, we write
    \[
    \mathcal{R}(M,g) := \{h \in \mathcal{R}(M) : h = \psi g, \psi \in C^\infty(M;(0,\infty)) \}
    \]
for the set of conformal perturbations of $g$.
We equip the set $\mathcal{R}(M,g)$ with the subspace topology inherited from $\mathcal{R}(M)$.

In this paper, we refer to a property  of $M$  being `generic' if the property holds for a residual set  of Riemannian metrics in $\mathcal{R}(M)$. 
 
We prove the following theorem.

\begin{namedthm*}{Main Theorem}\label{T:MainThm}
    Let $M$ be a connected compact smooth manifold of dimension $d \geq 2$ without boundary. 
    For a Riemannian metric $g$ on $M$, consider the following two properties:
    \begin{enumerate}
        \item the non-zero spectrum of $\Delta_g$ is linearly independent over $\mathbb{Q}$, and even stronger, all rational powers of non-zero eigenvalues
    \[
    \bigcup_{q \in \mathbb{Q}} \{\lambda_i^q : i \geq 1\}\subseteq \mathbb{R},
    \]
    are collectively linearly independent over $\mathbb{Q}$;
        \item the pointwise products of eigenfunctions $\{e^\alpha: \alpha\in\mathbb N^{(\infty)}\}$ are linearly independent.
    \end{enumerate} 
    
    The above two properties are generic, i.e.\ the set of Riemannian metrics $g \in \mathcal{R}(M)$ with the above two properties is residual.
    Moreover, for any Riemannian metric $g$ on $M$, the set of metrics in $\mathcal{R}(M,g)$ which have the above two properties is residual in $\mathcal{R}(M,g)$.
\end{namedthm*}

Most directly, this theorem is an extension of the result by Bando and Urakawa~\cite{BandoUrakawa1983} which provides that the metrics for which the spectrum of the Laplace--Beltrami operator is simple, is residual in both $\mathcal{R}(M)$ and $\mathcal{R}(M,g)$. 
Before outlining the strategy of the proof of our \autoref{T:MainThm}, let us make some remarks and point out a few consequences.

\begin{rem}\label{R:Strongest}
    Our method at its most general shows that for any  connected closed smooth  Riemannian manifold $(M,g)$  of dimension at least $2$ , and countable family of $C^1$-functions $\{F_k:\R^{N_k}\to \C\}_{k\in \N}$ with $\nabla F_k$ nowhere zero, there is a residual set of metrics $h \in \mathcal{R}(M,g)$ such that for all $k\in \N$, we have
    \[
    F_k(\lambda_1(h),\ldots, \lambda_{N_k}(h)) \neq 0.
    \]
    Likewise for any countable family of non-zero real-analytic functions $\{F_k:\R^{N_k}\to \C\}_{k\in \N}$.
\end{rem}

An important observation is that while the \autoref{T:MainThm} proves that generically the set 
\[
EE_N : = \mathrm{span} \{ e_i e_j : 0 \leq i \leq j \leq N-1\} \subseteq L_2(M),
\]
has dimension $\frac{N(N+1)}{2}$, a result by Lu, Sogge and Steinerberger~\cite{LuSogge2019} shows that $EE_N$ can be approximated well by a space with much smaller dimension.
We furthermore remark that with different methods, Privat and Sigalotti~\cite{PrivatSigalotti2010} had already proven that on generic domains, the Laplace operator with Dirichlet boundary conditions has linearly independent squares of eigenfunctions $\{e_i^2\}_{i \in \N}$, and spectrum linearly independent over $\mathbb{Q}$. 

A recent result which likewise deserves mention is a preprint by Greilhuber and Kepplinger~\cite{GreilhuberKepplinger2023} which is, in a way, complementary to our results. 
Where we limit ourselves to metrics with simple  Laplace--Beltrami  spectrum, they study the local neighbourhood of Riemannian metrics with higher multiplicity eigenvalues. 
This is in line with a research programme initiated by Arnol'd~\cite{Arnold1972,CdV1988}.

The following corollary of our \autoref{T:MainThm} might be relevant to noncommutative geometry:

\begin{cor}\label{C:Trunc}
    Let $(M,g)$ be a  closed smooth  Riemannian manifold  of dimension at least $2$.
    Represent the unital C$^*$-algebra $C(M;\mathbb{C})$ of complex-valued continuous functions on $M$ as bounded operators on $L_2(M,g)$. 
    Assume that $g$ belongs to the residual set of our \autoref{T:MainThm}.
    For $\Lambda \geq 0$, denote by $P_\Lambda = \chi_{[0,\Lambda]}(\sqrt{\Delta_g})$ the spectral projection of the square root of the Laplace--Beltrami operator onto the eigenspaces of eigenvalues at most $\Lambda$. 
    Then,
    \[
    P_\Lambda C(M)P_\Lambda \simeq \mathrm{Sym}_{N(\Lambda)}(\C),
    \]
    where $\mathrm{Sym}_{N(\Lambda)}(\C)$ denotes the set of matrices equal to their transpose, and $N(\Lambda)$ is the rank of $P_\Lambda$. 
\end{cor}
\begin{proof}
    Assume without loss of generality that the eigenfunctions $\{e_i\}_{i \in \mathbb{N}}$ are real-valued. The matrix elements of a multiplication operator $M_f \in \mathcal{B}(L_2(M,g))$, $f \in C(M;\C)$, are given by
    \[
        \langle e_n, M_f e_m\rangle = \langle e_n e_m, f\rangle = \langle e_m, M_f e_n\rangle.
    \]
    Since, by \autoref{T:MainThm}, we have that $\{e_i e_j : 0\leq  i \leq j \leq N(\Lambda)-1\}$ are linearly independent, the claim follows. 
\end{proof}

In recent years, a paradigm  popularised  by Connes and Van Suijlekom in noncommutative geometry on spectral truncations~\cite{ConnesvSuijlekom2021}
has attracted a lot of attention; from early investigations~\cite{DAndreaLizzi2014, GlaserStern2020} to the  more  recent~\cite{GlaserStern2021, DAndreaLandi2022, GaudillotSuijlekom2025, Rieffel2023, LeimbachvSuijlekom2024, Suijlekom2024}.
The core idea is to take a spectral triple~\cite{Connes1994} $(\Ac, \Hc, D)$ -- consisting of a  unital  $^*$-algebra of bounded operators on Hilbert space $\Hc$ and an  (essentially)  self-adjoint operator $D$ with compact resolvent,  which boundedly commutes with $\mathcal{A}$ 
-- and truncating this with a  spectral  projection $P = \chi(D)$ to $(P\Ac P, P\Hc, PDP)$. 
 The compression $P\Ac P$ generally only retains the structure of an operator system (i.e.\ a unital self-adjoint space of bounded operators on Hilbert space). 
However, various tools in noncommutative geometry are currently being adapted to operator systems, for a large part in order to investigate whether and how the spectral truncation $(P\Ac P, P\Hc, PDP)$ recovers information about the original spectral triple $(\Ac, \Hc, D)$.

\autoref{C:Trunc} shows that, for a generic connected closed smooth Riemannian manifold $M$, the compression $(PC(M)P, PL_2(M), P\sqrt{\Delta_g})$, $P = \chi_{[0,\Lambda]}(\sqrt{\Delta_g})$, does not retain any information about the spectral triple $(C(M),L_2(M,g),\sqrt{\Delta_g})$, except for the eigenvalues of modulus at most $\Lambda$.
This is in stark contrast to the highly symmetric situations investigated in \cite{Suijlekom2021, Rieffel2023, LeimbachvSuijlekom2024, GaudillotSuijlekom2025, Lei24, Suijlekom2024}, among others.
The expert will see quickly that \autoref{C:Trunc} implies that the spectral metric space $(C(M),L_2(M,g),\sqrt{\Delta_g})$ is not approximated by its spectral truncations in (quantum) Gromov--Hausdorff distance: 
in fact, any diagonal matrix in the compression $P_\Lambda C(M) P_\Lambda$ commutes with $P_\Lambda D P_\Lambda$, so if the rank of $P_\lambda$ is at least $2$, then the semi-norm $\lVert [P_\Lambda D P_\Lambda, \cdot ] \rVert$ on the compression will have non-trivial kernel, rendering the diameter of the state space of $P_\Lambda C(M) P_\Lambda$ with the Connes distance infinite.
That is to say, the pair $(P_\Lambda C(M) P_\Lambda, \lVert [P_\Lambda D P_\Lambda, \cdot ] \rVert)$ is not a compact quantum metric space in the sense of Rieffel (despite being finite-dimensional).

We admit that the spectral triple $(C(M),L_2(M,g),\sqrt{\Delta_g})$ is not the one which one would \emph{canonically} associate to a closed manifold:
instead the usual setting which noncommutative geometry tries to build upon is that of closed spin manifolds.
We will investigate the analogous genericity properties of Dirac spectra, as well as the consequences for the spectral truncation program, in a subsequent preprint. Finally, it also remains to be formally proven that there exist non-isometric isospectral manifolds inside the residual set of the \autoref{T:MainThm}.

\subsection*{Outline of the proof}
We base our methods on~\cite{BleeckerWilson1980, BandoUrakawa1983}, using conformal perturbations of the metric. 
First one readily notices that linear independence of the Laplace spectrum over $\mathbb{Q}$ is equivalent to linear independence over $\mathbb{Z}$.
For all $N, K, L \in \mathbb{N}$ with $K \geq L$, 
consider
the following two properties:
\begin{enumerate}
        \item $\lambda_1(g), \ldots, \lambda_N(g)$ are linearly independent over $\mathbb{Z}\cap[-K,K]$;
        \item the products $\{e^\alpha : \alpha \in \mathbb{N}^N \, |\alpha|\leq L\}$ are linearly independent.
\end{enumerate}
The intersection over all $N,K,L\in\mathbb{N}$, $K\geq L$, of the respective sets of Riemannian metrics with the above two properties is the set of metrics claimed to be residual in $\mathcal{R}(M)$ in our \autoref{T:MainThm} (one can similarly obtain a residual set for which rational powers of the eigenvalues are linearly independent over $\mathbb{Q}$ too).
Thus it remains to show that the above two properties hold for an open and dense set of Riemannian metrics in both $\mathcal{R}(M)$ and $\mathcal{R}(M,g)$.

As for openness of the first condition, note that it can be expressed as
    \[
    F_{k_1,\ldots,k_N}(g) := \sum_{i=1}^N k_i \lambda_i(g) \neq 0,
    \]
for all $k_i \in \mathbb{Z}\cap [-K,K]$. Since Bando--Urakawa \cite{BandoUrakawa1983} showed that the eigenvalues depend continuously on the Riemannian metric, this is therefore an open condition. For the products of eigenfunctions, around every $g\in \mathcal{R}(M)$ in the open set for which $\lambda_0 < \cdots < \lambda_N$, there exists an open neighbourhood for which the first $N$ eigenfunctions can be chosen to depend continuously on the metric~\cite[Lemma~3.2]{GreilhuberKepplinger2023}. Hence, via a Gram matrix argument, the combination of (1) and (2) is an open condition.

The main technical arguments happen in Section~\ref{S:MainProof}, where we show density of the above two conditions. In brief, we can summarise the strategy for dimension $d=2$ as follows. Along a conformal perturbation $g_t = (1+t\psi)g$, the eigenfunctions and eigenvalues depend real-analytically on $t$ (see Section~\ref{SS:Perturbation} below). For $d=2$, we have (Lemma~\ref{L:EigenvaluePerturbation})
\[
\dot{\lambda_i} =  -\lambda_i  \big\langle \psi, e_i^2\big\rangle_0.
\]
Hence, if we start with a metric $g$ for which the squares $\{e_i^2 : 0\leq i \leq N\}$ are linearly independent, we can find a function $\psi \in C^\infty(M)$ such that for $g_t = (1+t\psi)g$ condition~(1) is satisfied for all $0\neq t \in (-\varepsilon, \varepsilon)$ This gives density of condition (1), if the set of metrics for which the squares are linearly independent is dense. This is done in two steps:
\begin{enumerate}
    \item We reduce to the metrics $g \in \mathcal{R}(M)$ for which along \textit{each} perturbation $g_t = (1+t\psi)g$ and each linear dependency $\sum_{i=0}^N C_i e_i^2 \equiv 0$, we can provide an extension real-analytically in $t$ for a short time-interval $t\in(-\varepsilon, \varepsilon)$
    \[
    \sum_{i=0}^N C_i(t) e_i(t)^2 \equiv 0.
    \]
    \item We show that if $g \in \mathcal{R}(M)$ is such as in (1), then 
    \[
    \sum_{i=0}^N C_i e_i^2 \equiv 0 \qquad \Rightarrow \qquad \sum_{i=0}^N \lambda_i C_i e_i^2\equiv 0.
    \]
    If $\lambda_0 < \lambda_1< \cdots< \lambda_N$, this forces $C_i=0$ for all $i$.
\end{enumerate}
Step (2) is achieved by taking the derivative in $t$ of the expression $\sum_{i=0}^N C_i(t) e_i(t)^2 \equiv 0$ and observing the implications for this derivative to be zero for each $\psi \in C^\infty(M)$. Finally, linear independence of the higher products $e^\alpha$ is obtained by the same strategy, though this now requires that sums $\sum_{j=1}^N \alpha_j \lambda_j$ take unique values. However, that follows from the spectrum being linearly independent over $\mathbb{Z}\cap[-K,K]$ for $K$ large enough.

\subsection*{Standing assumptions and notation}
Throughout this paper, we will exclusively be working with manifolds which are connected, smooth and without boundary.
We will usually additionally assume $M$ to be compact and of dimension $d \geq 2$.
Usually $M$ will be equipped with a Riemannian metric $g$.
Whenever convenient we will use the Einstein summation convention.
We write $C^\infty(M)$ for the set of all \emph{real}-valued smooth functions on $M$.
For smooth vector fields $X,Y$ on $M$, we denote by $X \cdot Y$ the smooth function given by $g(X,Y) = g_{ij} X^i Y^j$. 
Recall that a smooth vector field $X$ on $M$ acts on a smooth function by $X(f) := X \cdot \nabla_gf$, where $\nabla_gf$ is the gradient vector field of $f$.

\section*{Acknowledgments}

The second and third named authors benefited from support by the Max Planck Institute for Mathematics in Bonn. Furthermore, we acknowledge helpful discussions with Edward McDonald.

\subsection*{AI usage}

In writing this paper, some conversations with Claude model Opus 4.8 and 5 were helpful in scoping the literature on genericity results and perturbation techniques. In terms of the current proof strategy, the idea of analysing the perturbation formula~\eqref{eq:Absurdity2} as a pseudodifferential operator must be credited to AI. It furthermore assisted with computing some of the symbols appearing in Section~\ref{S:MainProof}. Apart from proofreading and minor notational improvements with Claude model Opus 5.5, no AI was used in the actual writing of the text and proofs in this paper.

\section{Preliminaries}

\subsection{Hilbert space lemmas}

Throughout this subsection, let $\Hc$ be a complex Hilbert space with inner product $\langle\cdot,\cdot\rangle$ anti-linear in the first component.

\begin{lem}\label{L:LinIndep}
    Given $v_1,\ldots, v_n \in \Hc$ linearly independent, and $c_1,\ldots,c_n \in \C$, we can find a unique $v\in \mathrm{span}\{v_1,\ldots,v_n\}$ such that $\langle v,v_i\rangle=c_i$ for $i=1,\ldots,n$.
\end{lem}
\begin{proof}
    It suffices to prove the lemma for $\Hc=\C^n$. Let $\{e_1,\ldots,e_n\}$ be the standard basis of $\C^n$. Then the conditions $\langle v,v_i\rangle =\sum_{j=1}^n \langle v, e_j\rangle \langle e_j, v_i\rangle= c_i$ form a system of $n$ linear equations for the variables $\langle v,e_j\rangle$. This has a unique solution found by Gaussian elimination due to the linear independence of the vectors $v_1,\ldots, v_n$.
\end{proof}

\begin{lem}\label{L:GramMatrix}
    Consider a collection of real-analytic maps into a Hilbert space
    \begin{align*}
    v_i: (-\varepsilon, \varepsilon) \to \Hc, \quad i\in \{1,\ldots, N\}.
    \end{align*} Let $G(t) \in M_N(\C)$ be the Gram matrix corresponding to these vectors, i.e. $G(t) = (\langle v_i(t), v_j(t)\rangle)_{i,j=1}^N$. Then after possibly shrinking the domain interval, we must have one of the following:
    \begin{enumerate}
        \item $\dim \ker G(t)$ is constant;
        \item $\dim \ker G(t) = k < \dim \ker G(0)$ for $0\neq t\in (-\varepsilon, \varepsilon)$.
    \end{enumerate}
    In the first case, 
    for any vector $C \in \C^N$ with
    \[
    \sum_{i=1}^N C_i v_i=0,
    \]
    there is
    a real-analytic map $t \mapsto C(t) \in \C^N$ such that 
    $C(0)=C$ and
    \[
    \sum_{i=1}^N C_{i}(t) v_i(t) = 0, \quad t \in(-\varepsilon,\varepsilon).
    \]
\end{lem}
\begin{proof}
    The Gram matrix $G(t)$ is positive semi-definite and with real-analytic entries.
    By a theorem due to Rellich~\cite[Satz~I]{Rellich1937a}, we can find real-analytic maps $\lambda_j(t) \in \C$ and $e_i(t) \in \C^N$, $1\leq i \leq N$, such that
    \begin{align*}
        G(t) e_i(t) &= \lambda_i(t)e_i(t);\\
        \langle e_i(t), e_j(t)\rangle &= \delta_{ij}. 
    \end{align*}
    Each real-analytic map $\lambda_i : (-\varepsilon, \varepsilon)\to \R$ can only have finitely many zeroes, unless it is identically zero. By shrinking the domain interval, we can assume without loss of generality that the only possible zeroes are at $t=0$. 
    This means precisely that either $\dim \ker G(t) = \dim \ker G(0)$ for each $t$ in the domain, or $\dim \ker G(t) = k < \dim \ker G(0)$ for $t\neq 0$. 

    Assume now that $\dim \ker G(t)$ remains constant, so that that the real-analytic maps $\lambda_j(t)$ with $\lambda_j(0)=0$ are constant and the corresponding eigenfunctions $e_j(t)$ remain in the kernel of $G(t)$. Suppose we have $C\in \C^N$ with
    \[
    \sum_{j=1}^N C_j v_j = 0.
    \]
    By construction, this means that $G(0)C = 0$. Decompose $C$ into the eigenvectors $e_j(0)\in \C^N$ spanning the kernel of $G(0)$,
    \[
    C = \sum_{\lambda_j=0} \tilde{C}_j e_j(0),
    \]
    so that $C$ also extends analytically to a function $C(t) \in \C^N$ given by 
    \[
    C(t) := \sum_{\lambda_j=0} \tilde{C}_j e_j(t).
    \]
    Write now $C_j(t)$ for the $j$-th component of $C(t) \in \C^N$, so that $C_j(0) = C_j$. We have $G(t)C(t) = 0$, and writing out the matrix multiplication means precisely that
    \begin{align*}
        0
        =\sum_{j=1}^N \langle v_i(t),v_j(t)\rangle C_j(t)
        = \langle v_i(t), \sum_{j=1}^N C_j(t)v_j(t) \rangle. 
    \end{align*}
    Hence, for all $t \in (-\varepsilon, \varepsilon)$, we have that 
    \[
    \sum_{j=1}^N C_j(t)v_j(t) \in \mathrm{span} \{v_i(t) : 1\leq i \leq N\} \cap \big(\mathrm{span} \{v_i(t) : 1\leq i \leq N\}\big)^{\perp} = \{0\}.\qedhere
    \]
\end{proof}

\subsection{Perturbation theory}\label{SS:Perturbation}
First of all, without loss of generality we will always assume that the eigenfunctions $\{e_i\}_{i \in \mathbb{N}}$ of the Laplace--Beltrami operator are real-valued. 

We will make extensive use of the following two propositions. We start with a result by Besson~\cite{Besson1989}. Note that Besson proved his result for $C^k$-metrics, a proof for the version for smooth metrics can be found in~\cite{GreilhuberKepplinger2023}. 
We will only need a special case for simple eigenvalues. 
 \begin{prop}\label{P:BessonSubmersion}
     Let $(M,g)$ be a closed Riemannian manifold of dimension at least $2$, and let $\lambda$ be a simple eigenvalue of $\Delta_g$. Then there exists an open neighbourhood $U$ of $g$ in $\mathcal{R}(M)$ and a $C^1$ submersion $\pi: U \to \R$ and $\varepsilon>0$ such that 
     \[
     \sigma(\Delta_{h}) \cap (\lambda-\varepsilon, \lambda+\varepsilon) = \pi(h),
     \]
     for all $h \in U$.
 \end{prop}

Next, we provide a proposition due to Berger~\cite{Berger1973}. Again, we only provide the version for simple eigenvalues.
\begin{prop}\label{P:AnalyticDependence}
    Let $(M,g_0)$ be a compact Riemannian manifold. For $\psi \in C^\infty(M; \mathbb{R})$, put $g_t = (1+t\psi)g_0$, and let $t_0 > 0$ such that $g_t$ is a smooth Riemannian metric for $|t|< t_0$. For every simple eigenvalue $\lambda$ of $\Delta_{g_0}$, there exist $\lambda(t) \in \R$ and $e(t)\in C^\infty(M; \mathbb{R})$ such that for $|t|<t_0$:
    \begin{itemize}
    \item $\lambda(0) = \lambda$;    
        \item $\lambda(t)$ and $e(t)$ depend real analytically on $t$;
        \item $\Delta_{g_t}e(t) = \lambda(t) e(t)$;
        \item $\|e(t)\|^2_{L_2(M, g_t)} =1$.
    \end{itemize}
    Furthermore, after possibly shrinking the domain to $t\in (-\varepsilon,\varepsilon)$, 
    \begin{itemize}
        \item there exists $\delta>0$ such that $\#(\sigma(\Delta_t) \cap (\lambda - \delta, \lambda +\delta)) = 1$. In other words, the function $\lambda(t)$ does not cross other eigenvalues.
    \end{itemize}
\end{prop}
\begin{proof}
     The first part of the lemma is~\cite[Lemma~3.15]{Berger1973}, for an explicit proof in the current setting see~\cite[Appendix]{BandoUrakawa1983}. The second part concerning $\lambda(t)$ not crossing other eigenvalues follows from Proposition~\ref{P:BessonSubmersion}.
\end{proof}

We now provide three well-known results on the derivatives of these analytic eigenvalues and eigenfunctions.
The following lemma can be found in e.g.~\cite[Lemma~4.2]{BandoUrakawa1983}, but we provide a self-contained proof.

\begin{lem}\label{L:ConformalLaplacian}
    Let $(M,g_0)$ be a compact $d$-dimensional Riemannian manifold, and let $\psi \in C^\infty(M,\mathbb{R})$. 
    Let $\varepsilon > 0$ be small enough such that $t\psi > -1$, for all $t \in (-\varepsilon,\varepsilon)$.
    Consider the conformal metric $g_t= (1+t\psi) g_0$, for all $t\in(-\varepsilon,\varepsilon)$.
    Then the conformal Laplace--Beltrami operator is given by
    \begin{align}
    \label{eqn:ConformalLaplacian} \Delta_t f &= \frac{1}{1+t\psi}\Delta_0 f  -\frac{d-2}{2}  \frac{t}{(1+t\psi)^2} \nabla_{g_0} \psi \cdot \nabla_{g_0} f, \quad f \in C^\infty(M).
    \end{align}
    We record furthermore that
    \[
    \dot{\Delta}_0 f := \frac{d}{dt}\Big|_{t=0} \Delta_t(f) = -\psi \Delta_0 f - \frac{d-2}{2}\nabla_{g_0}\psi \cdot \nabla_{g_0}f.
    \]
\end{lem}

\begin{proof}
    For $t\in(-\varepsilon,\varepsilon)$, set $\phi_t := 1+t\psi \in \rmC^\infty(M,\R_+)$.
    Denote $|g_t| := \det g_t$, and recall $|g_t| = \phi_t^d |g_0|$.
    Then for $f \in \rmC^\infty(M)$, working in local coordinates, and using the product rule for $\partial_i$, we have:
    \begin{align*}
        \Delta_tf
        &= -\frac{1}{|g_t|^{\frac{1}{2}}} \partial_i \Big( |g_t|^{\frac{1}{2}} g_t^{ij} \partial_j f \Big) \\
        &= -\frac{1}{\phi_t^{\frac{d}{2}} |g_0|^{\frac{1}{2}}} \partial_i \Big( (\phi_t^{\frac{d}{2}} |g_0|^{\frac{1}{2}} \phi_t^{-1} g_0^{ij} \partial_j f \Big) \\
        &= -\frac{1}{\phi_t^{\frac{d}{2}}  |g_0|^{\frac{1}{2}} } \Big( \frac{d-2}{2} \phi_t^{\frac{d}{2}-2} (\partial_i\phi_t) |g_0|^{\frac{1}{2}} g_0^{ij} \partial_j f + \phi_t^{\frac{d}{2}-1} \partial_i \big( |g_0|^{\frac{1}{2}} g_0^{ij} \partial_j f \big) \Big) \\
        &= -\frac{d-2}{2} \frac{1}{\phi_t^2} g_0\big(\nabla_{g_0}\phi_t,\nabla_{g_0}f\big) + \frac{1}{\phi_t} \Delta_0f.
    \end{align*}
\end{proof}

\begin{lem}\label{L:EigenfunctionDerivative}
    Let $(M,g_0)$ be a compact $d$-dimensional Riemannian manifold, and let $\psi \in C^\infty(M)$. 
    Consider the metric $g_t= (1+t\psi) g_0$ for $t$ small enough. 
    Assume that the eigenvalue $\lambda_i$ is simple. 
    Then, writing $\dot{e}_i := \frac{d}{dt}\big|_{t=0} e_i(t)  \in C^\infty(M;\R)$, we have that
    \begin{align*}
    \dot{e}_i &= -\frac{d}{4}\Big(\int_M e_i^2 \psi\, d\nu_0 \Big) e_i + \sum_{k\neq i} \frac{\lambda_i}{\lambda_k-\lambda_i} \Big(\int_M e_ke_i  \psi\, d\nu_0\Big) e_k\\
    &\quad + \Big(\frac{d}{2}-1\Big) \sum_{k\neq i} \frac
    {1}{\lambda_k - \lambda_i}\bigg( \int_M e_k \nabla_g\psi \cdot \nabla_ge_i \, d\nu_0\bigg)e_k\\
    &= -\frac{d}{4}\langle e_i^2,\psi\rangle_0\,  e_i + \lambda_i f_i(\Delta_0)(e_i\psi) + \frac{d-2}{2} f_i(\Delta_0) \Big(\nabla_g\psi \cdot \nabla_ge_i \Big),
    \end{align*}
    where 
    \[
    f_i(x) = \begin{cases}
        0 & x=\lambda_i,\\
        \frac{1}{x-\lambda_i} & x\neq \lambda_i.
    \end{cases}
    \]
\end{lem}
\begin{proof}
The derivatives $\dot{e}_i$ are smooth real-valued functions by Proposition~\ref{P:AnalyticDependence}, and hence they are in $L_2(M, g_0)$. It therefore suffices to compute the inner products $\langle e_k, \dot{e}_i\rangle_0$.

    Recall 
    that $d\nu_t = (1+t\psi)^{\frac{d}{2}}d\nu_0$.
    The condition  from \autoref{P:AnalyticDependence}  that the eigenfunctions $e_i(t)$ stay normalised, therefore means that
    \[
    \int_M e_i(t)^2 (1+t\psi)^{\frac{d}{2}}\, d\nu_0 = 1.
    \] 
    Taking derivatives in $t$ on both sides and evaluating at $t=0$, we obtain
    \begin{align*}
        2\int_M e_i \dot{e}_i \, d\nu_0 + \frac{d}{2}\int_M e_i^2 \psi\, d\nu_0
        =0,
    \end{align*}
    which gives that
    \[
    \langle e_i, \dot{e}_i \rangle_{g_0} = 
    -\frac{d}{4} \int_M e_i^2 \psi\, d\nu_0.
    \]

Next, by (\ref{eqn:ConformalLaplacian}), we have that 
\begin{align*}
    \Delta_t e_i(t) &= \lambda_i(t)e_i(t)\\
    \iff \Delta_0e_i(t) &=  \frac{d-2}{2}\frac{t}{1+t\psi} \nabla_g \psi \cdot \nabla_g e_i(t) + (1+t\psi) \lambda_i(t)e_i(t).
\end{align*}
Now take the inner product (with respect to $d\nu_0$) with $e_k$ (not depending on $t$) for $k\neq i$, and we obtain
\begin{align*}
    \langle e_k, \Delta_0e_i(t)\rangle_0 &= \frac{d-2}{2}\langle e_k,\frac{t}{1+t\psi} \nabla_g \psi \cdot \nabla_g e_i(t)\rangle_0+\lambda_i(t) \langle e_k,(1+t\psi) e_i(t)\rangle_0.
\end{align*}
Observe that $\langle e_k, \Delta_0e_i(t)\rangle_0=\lambda_k \langle e_k, e_i(t)\rangle_0$ since $\Delta_0$ is symmetric.  
Now again take derivatives in $t$ and evaluate at $t=0$, and we obtain
\begin{align*}
    \lambda_k\langle e_k, \dot{e}_i\rangle_0 &=  \frac{d-2}{2}\langle e_k, \nabla_g \psi \cdot \nabla_g e_i\rangle_0+\dot{\lambda}_i \underbrace{\langle e_k, e_i\rangle_0}_{=0} + \lambda_i \langle e_k, \psi e_i\rangle_0 + \lambda_i \langle e_k, \dot{e}_i\rangle_0.
\end{align*}
Hence indeed,
\[
\langle e_k, \dot{e}_i\rangle_0 = \frac{\lambda_i}{\lambda_k-\lambda_i} \int_M e_k e_i  \psi\, d\nu_0 + \frac{d-2}{2} \frac{1}{\lambda_k-\lambda_i}\int_M e_k \nabla_g \psi \cdot \nabla_ge_i\,d\nu_0.
\]
    Finally, note that
    \begin{align*}
        \sum_{k\neq i} \frac{1}{\lambda_k -\lambda_i} \Big( \int_M e_i e_k \psi\, d\nu_0 \Big)e_k &=  \sum_{k\neq i} \frac{1}{\lambda_k -\lambda_i} \langle e_k, e_i \psi \rangle_0\, e_k\\
        &= \sum_{k \in \N} \langle f_i(\Delta_0)e_k, e_i \psi \rangle _0\, e_k \\
        &= f_i(\Delta_0) (e_i\psi), 
    \end{align*}
    where  $f_i(\Delta_0) = (\Delta-\lambda_i)^{-1}(1-P_{\{\lambda_i\}})$, i.e.
      \begin{align}\label{eqn:Resolvent-i}
    f_i(x) = \begin{cases}
        0 & x=\lambda_i,\\
        \frac{1}{x-\lambda_i} & x\neq \lambda_i.
    \end{cases}
    \end{align}
    Likewise for the other term.
\end{proof}

We now briefly remark that $f_i(\Delta_0)$ appears in~\cite[Section~2]{Zelditch1990} in the notation $G_\lambda^\perp$, work by Zelditch which shows the existence of isospectral manifolds with simple spectrum. 

We furthermore provide the derivative of the eigenvalues.

\begin{lem}\label{L:EigenvaluePerturbation}
    Let $(M,g_0)$ be of dimension $d\geq 2$. 
    Consider the metric $g_t= (1+t\psi) g_0$ for $t$ small enough. Suppose that $\lambda_i$ is an eigenvalue of $\Delta_0$, and consider an analytic family $\lambda_i(t)$, $e_i(t)$ as in Proposition~\ref{P:AnalyticDependence}, such that $\Delta_t e_i(t) = \lambda_i(t)e_i(t)$. Then,
    \[
    \dot{\lambda_i} =    \big\langle \psi, \Big(\frac{2-d}{4}\Delta_0-\lambda_i\Big)e_i^2\big\rangle_0.
    \]
\end{lem}
\begin{proof}
    Differentiate at $t=0$ the identity
    \[
    \langle \Delta_t e_i(t), e_i(t)\rangle_0 = \lambda_i(t) \langle e_i(t), e_i(t) \rangle_0, 
    \]
    using Lemma~\ref{L:ConformalLaplacian}, and we obtain
    \begin{align*}
        -\lambda_i\langle \psi e_i, e_i\rangle_0 -\frac{d-2}{2} \langle \nabla_g \psi \cdot \nabla_ge_i, e_i\rangle_0+2\lambda_i\langle e_i, \dot{e}_i\rangle_0 &= \dot{\lambda_i}+2\lambda_i\langle e_i, \dot{e}_i\rangle_0.
    \end{align*}
    Now note that, just as in~\cite[Proposition~4.6]{BleeckerWilson1980},
    \begin{align*}
    2\langle \nabla_g \psi \cdot \nabla_ge_i, e_i\rangle_0 &= \langle (\Delta_0 \psi) e_i, e_i\rangle_0 + \langle \psi \Delta_0 e_i, e_i\rangle_0 - \langle \Delta_0(\psi e_i),e_i\rangle_0\\
    &= \langle \Delta_0 \psi , e_i^2\rangle_0 + \lambda_i \langle \psi e_i, e_i\rangle_0   -\lambda_i \langle \psi e_i, e_i\rangle_0\\
    &= \langle \Delta_0\psi, e_i^2\rangle_0.
    \end{align*}
    Assembling these computations provides the assertion.
\end{proof}

Finally, we end with a lemma which will become useful later.

\begin{lem}\label{L:CosphereIntegral}
    For $f \in C^\infty(M)$, we have
    \begin{align*}
        \int_{|\xi|_g=1} \nabla^2f(\xi,\xi)\, d\xi &= \frac{1}{d} \tr_g \nabla^2 f = -\frac{1}{d}\Delta_g f;\\
        \int_{|\xi|_g=1} (\xi \cdot \nabla_g f)^2\, d\xi &= \frac{1}{d} |\nabla_g f|^2.
    \end{align*}
\end{lem}
\begin{proof}
    Both formulas follow from doing the computation in geodesic normal coordinates, combined with the fact that for the normalised volume form on the sphere $\mathbb{S}^{d-1}$, we have~\cite{Folland2001}
    \[
    \int_{\mathbb{S}^{d-1}} \xi_i \xi_j \, d\xi = \frac{1}{d} \delta_{ij}. 
    \]
\end{proof}

\section{Main proof}\label{S:MainProof}
From the computations in Lemma~\ref{L:EigenvaluePerturbation}, along a conformal perturbation of the metric we have 
    \[
    \dot{\lambda_i} =    \big\langle \psi, \Big(\frac{2-d}{4}\Delta_0-\lambda_i\Big)e_i^2\big\rangle.
    \]
To obtain complete control over these derivatives, we will prove that for a dense set of metrics $g$, the functions $\{\big(\frac{2-d}{4}\Delta_g-\lambda_i\big)e_i^2 : 1\leq i \leq N\}$ are linearly independent. This will be derived from density of the set of metrics for which the squares $\{e_i^2\}_{i \in \mathbb{N}}$ are linearly independent. There it will turn out to be crucial that the spectrum of the Laplace--Beltrami operator is generically simple. Finally, to prove the generic linear independence of the products $\{e^\alpha : \alpha \in \mathbb{N}^{(\infty)}\}$, it will be crucial that the sums $\sum_{j=1}^N \alpha_j \lambda_{j}$ of eigenvalues are generically non-degenerate. The proof strategy in this section is therefore necessarily a chain of implications:
\begin{enumerate}
    \item density of the metrics for which $\{e_i(g)^2: 0\leq i \leq N\}$ is linearly independent;
    \item density of the  metrics for which $\{ \big(\frac{2-d}{4}\Delta_g-\lambda_i(g)\big)e_i(g)^2: 1 \leq i \leq N\}$ is linearly independent;
    \item density of the  metrics for which sums $\{\sum_{j=1}^N \alpha_j\lambda_{j}(g) : \alpha \in \mathbb{N}^N,\, |\alpha|\leq L\}$, are non-degenerate;
    \item density of the  metrics for which $\{e(g)^\alpha : \alpha \in \mathbb{N}^N,\, |\alpha|\leq L\}$ is linearly independent.
\end{enumerate}

\subsection{Linear independence of squares}\label{SS:Squares}
We will first aim to obtain linear independence of the squares of the eigenfunctions. 
Later we will use the same methods, so this section also serves as a stepping stone to help digest later proofs.

\begin{thm}\label{T:SquareIndependence}
    Let $M$ be a compact manifold. 
    For each $N\in \N$, the set of metrics $g$ in $\mathcal{R}(M)$ for which simultaneously 
     
    the following two properties hold, is open and dense in $\mathcal{R}(M)$:
    
    \begin{enumerate}
        \item $0=\lambda_0(g) < \lambda_1(g)<\cdots<\lambda_N(g) $;
        \item $\{e_i^2 : 0\leq i \leq N\}$ is a set of linearly independent functions.
    \end{enumerate}
\end{thm}

Recall that~\cite{BandoUrakawa1983} proves precisely that the set of metrics for which $0=\lambda_0(g) < \lambda_1(g)<\cdots<\lambda_N(g)$ is open and dense in $\mathcal{R}(M)$.

\begin{lem}\label{L:Hammer1}
    Let $M$ be a smooth compact manifold and let $N\in  \N$. 
    Then, the set of metrics for which simultaneously
     
    the following two properties hold, is dense in $\mathcal{R}(M)$:
    
    \begin{enumerate}
        \item $\lambda_0(g) < \lambda_1(g)<\cdots<\lambda_N(g) $;
        \item for any $\psi \in C^\infty(M)$ giving rise to metrics $g_t = (1+t\psi)g$ and for any $C_i \in \C$ for which
    \[
    \sum_{i=0}^N C_i e_i^2 \equiv 0,
    \]
    there 
    exists a (small) positive real number $\varepsilon$ and $\C$-valued analytic functions $(-\varepsilon,\varepsilon) \ni t \mapsto C_i(t)$ with $C_i(0)=C_i$, such that
    \[
    \sum_{i=0}^N C_i(t) e_i(t)^2 \equiv 0,
    \]
     for all $t \in (-\varepsilon,\varepsilon)$.
    \end{enumerate} 
\end{lem}
\begin{proof}
    By~\cite{BandoUrakawa1983}, the set of metrics for which $\lambda_0(g) < \lambda_1(g)<\cdots<\lambda_N(g)$ is dense in $\mathcal{R}(M)$. Hence, take any $g_0 \in \mathcal{R}(M)$ with $\lambda_0(g_0) < \cdots < \lambda_N(g_0)$. Along any conformal perturbation $(1+t\psi)g_0$, one can pick eigenfunctions $e_0(t)$, $1\leq i \leq N$ to depend analytically on $t$ due to Proposition~\ref{P:AnalyticDependence}.

    Now by Lemma~\ref{L:GramMatrix}, after possibly shrinking the domain interval, either $\dim \ker G(t) < \dim \ker G(0)$ for $0 \neq t \in (-\varepsilon, \varepsilon)$, or $\dim \ker G(t) = \dim \ker G(0)$ for all $t \in (-\varepsilon, \varepsilon)$. In the first case, replace $g_0$ by $g_\delta$ for some small $\delta$. Iterate until $\dim \ker G(t)$ is constant along any perturbation for some small time $t$. By this procedure, we find a metric $g$ which lies as close to $g_0$ as we like, and for which $\dim \ker G(t)$ is constant in any conformal perturbation direction for some small time (depending on $\psi$). Finally, Lemma~\ref{L:GramMatrix} provides that if $\dim \ker G(t)$ remains constant along a conformal perturbation, any relation $\sum_{i=0}^N C_i e_i^2 \equiv 0$ can be extended analytically to $\sum_{i=0}^N C_i(t) e_i(t)^2 \equiv 0$.
\end{proof}

We will now show that in the presence of simple spectrum, the real-analytic continuations  $\sum_{i=0}^N C_i(t) e_i(t)^2 \equiv 0$ can only exist along all conformal perturbations if $C_i=0$ for each $i$. 
Take the derivative at $t=0$ of the identity $\sum_{i=0}^N C_i(t) e_i(t)^2 \equiv 0$ appearing in Lemma~\ref{L:Hammer1}, as all terms are analytic in $t$:
\begin{equation*}\label{eq:Absurdity1}
\sum_{i=0}^N \dot{C}_i e_i^2 = -2\sum_{i=0}^N C_i e_i \dot{e}_i.
\end{equation*}
Using the expression for $\dot{e_i}$ in Lemma~\ref{L:EigenfunctionDerivative}, this identity expands as
\begin{equation}\label{eq:Absurdity2}\begin{split}
\sum_{i=0}^N \Big(\dot{C}_i - \frac{d}{2} C_i \langle e_i^2,\psi\rangle_0 \Big)  e_i^2 &= -2\sum_{i=0}^N \lambda_i C_i f_i(\Delta_0)(e_i \psi) e_i -(d-2)\sum_{i=0}^N C_i f_i(\Delta_0)(\nabla_g\psi \cdot \nabla_g e_i) e_i.   
\end{split}
\end{equation}

The left-hand side lies in the span of $\{e_i^2: 0\leq i \leq N\}$, which is finite-dimensional. The right-hand side can be regarded as a classical pseudodifferential operator applied to $\psi$, namely (note that the $i=0$ term vanishes, as $\lambda_0=0$ and $\nabla_g e_0=0$)
\begin{equation}\label{eq:PsDOT}
T:= -2\sum_{i=1}^N \lambda_i C_i M_{e_i} \circ f_i(\Delta_0) \circ M_{e_i} -(d-2)\sum_{i=1}^N C_i  M_{e_i}\circ f_i(\Delta_0) \circ \nabla_g e_i \in \Psi^{-1}_{cl}(M), 
\end{equation}
where $M_{e_i}$ is a multiplication operator, and $ \nabla_g e_i$ is regarded as a vector field on $M$ and hence a first-order differential operator. 

\begin{lem}
    If $T\in \Psi(M)$ has finite-rank, then $T \in \Psi^{-\infty}(M)$.
\end{lem}
\begin{proof}
    Since $T: C^\infty(M) \to C^\infty(M)$ has finite-rank, we obtain by that $T: \Hc^s\to C^\infty(M)$ has finite rank for all $s\in \R$, where $\Hc^s$ is the Bessel potential Sobolev space. This  implies $T\in \Psi^{-\infty}(M)$.
\end{proof}

Hence, equation~\eqref{eq:Absurdity2} can only hold for all $\psi \in C^\infty(M)$ if the asymptotic expansion of the symbol of the classical pseudodifferential operator $T$ into homogeneous parts
\[
\sigma(T) \sim \sum_{j\geq 0} \sigma_{-1-j},
\]
satisfies $\sigma_{-1-j}(x,\xi) = 0$ everywhere.

\begin{lem}\label{L:TPrincipalSymb}
    For $T \in \Psi^{-1}_{cl}(M)$ defined in equation~\eqref{eq:PsDOT}, we have for $x$ in the centre of geodesic normal coordinates,
    \begin{align*}
        \sigma_{-1}(T)(x,\xi) &= -i(d-2)\sum_{i=1}^N C_i e_i(x) \frac{\xi \cdot \nabla_ge_i}{|\xi|_g^2};\\
        \sigma_{-2}(T)(x,\xi)&= -2\sum_{i=1}^N \lambda_i C_i e_i(x)^2\frac{1}{|\xi|_g^2} +2(d-2)\sum_{i=1}^N C_i e_i(x)\frac{\nabla^2e_i(\xi,\xi)}{|\xi|_g^4}.
    \end{align*}
\end{lem}
\begin{proof}
    Recall that, for any smooth function $f \in C^\infty(M)$, the symbols of the multiplication operator $M_f$ (considered as a pseudodifferential operator of order $0$), and of the gradient vector field $\nabla_gf$ (considered as a pseudodifferential operator of order $1$, acting by $\nabla_gf(h)(x) := \nabla_gf(x) \cdot \nabla_gh(x)$), are respectively given by
    \begin{align*}
        \sigma(M_f)(x,\xi) &= f(x), \; \text{and} \\
        \sigma(\nabla_gf)(x,\xi) &= i\xi(\nabla_g f(x)) = ig^{kl} \xi_k \partial_{x_l} f(x),
    \end{align*}
    for all $x \in M$ and $\xi \in T^*_xM$.
    Here and in the sequel we are using the Einstein summation convention. Furthermore, recall that in local coordinates, we have
    \[
    \Delta_0 f(x) = - g^{ij}(x) \partial_i \partial_j f(x) + g^{ij}(x) \Gamma^{k}_{ij}(x) \partial_k f(x),
    \]
    and so in these coordinates, we have
    \begin{equation}\label{eq:HomPartsLaplace}
    \sigma(\Delta_0)(x,\xi) = g^{ij}(x)\xi_i\xi_j + ig^{ij}(x)\Gamma^{k}_{ij}(x)\xi_k = : p_2(x,\xi) + p_1(x,\xi).
    \end{equation}
    This expression is dependent on the coordinate chart.

    Now recall that $f_i(\Delta_0) := (1-P_{\{\lambda_i\}})(\Delta_0-\lambda_i)^{-1}$. We will obtain the expansion of the symbol of $f_i(\Delta_0)$ in a chart into homogeneous parts from that of $\Delta_0-\lambda_i$. We have that $f_i(\Delta_0) \circ (\Delta_0 - \lambda_i) =1-P_{\{\lambda_i\}}$, and so if in a chart
    \begin{align*}
        \sigma(\Delta_0-\lambda_i) &\sim p_{2} + p_1 -\lambda_i \\
        \sigma(f_i(\Delta_0))&\sim r_{-2} + r_{-3} + \cdots,
    \end{align*}
    then in the same chart,
    \begin{align*}
        \sigma(f_i(\Delta_0)\circ (\Delta_0-\lambda_i)) &\sim \sum_{|\alpha|=0}^\infty \frac{(-i)^{|\alpha|}}{\alpha!}\partial_\xi^\alpha \sigma(f_i(\Delta_0)) \partial_x^\alpha \sigma(\Delta_0-\lambda_i)\\
        &\sim (p_2 r_{-2}) + \Big(p_1 r_{-2} + p_2r_{-3}-i \sum_{|\alpha|=1} \partial_\xi^\alpha  r_{-2} \partial_x^\alpha  p_2  \Big)+ \cdots,
    \end{align*}
    where the grouping is respectively the order $0$ and order $-1$ homogeneous part of the symbol of $f_i(\Delta_0)\circ (\Delta_0-\lambda_i)$. 
    However, since $f_i(\Delta_0)\circ (\Delta_0-\lambda_i) = 1-P_{\{\lambda_i\}}$,  i.e.\ the difference of the identity and a smoothing operator, we deduce that 
    \begin{align*}
    1&= p_2 r_{-2};\\
    0&= p_1 r_{-2} + p_2r_{-3}-i \sum_{j=1}^d \partial_{\xi_j} r_{-2} \partial_{x_j} p_2 .
    \end{align*}
    Fix geodesic normal coordinates centred at $x$. Then from the explicit expressions for $p_1$ and $p_2$ in equation~\eqref{eq:HomPartsLaplace}, we obtain
    \begin{align*}
        r_{-2}(y,\xi) &= \frac{1}{p_2}= \frac{1}{|\xi|^2_{g(y)}}\\
        r_{-3}(y,\xi)&= \frac{1}{p_2}\Big(-p_1r_{-2}+i\sum_{j=1}^d \partial_{\xi_j} r_{-2}  \partial_{x_j} p_2 \Big)\\
        &= -\frac{ig^{ij}(y)\Gamma^{k}_{ij}(y)\xi_k}{|\xi|^4_{g(y)}}+\frac{i}{|\xi|_{g(y)}^2}\sum_{j=1}^d \partial_{\xi_j}  \Big(\frac{1}{ g^{mn}(y)\xi_m\xi_n}\Big) \partial_{x_j}  \Big(g^{kl}(y)\xi_k\xi_l \Big)\\
        &=-\frac{ig^{ij}(y)\Gamma^{k}_{ij}(y)\xi_k}{|\xi|^4_{g(y)}}-2i\sum_{j=1}^d g^{mj}(y)\xi_m  \frac{\partial_{x_j}g^{kl}(y)\xi_k\xi_l}{|\xi|_{g(y)}^6},
    \end{align*}
    where in the last step, we used that
    \begin{align*}
     \partial_{\xi_j} \Big(\frac{1}{ g^{mn}(y)\xi_m\xi_n}\Big) 
     &= - \frac{1}{(g^{mn}(y)\xi_m\xi_n)^2} 2 g^{mj}(y) \xi_m 
     =  -2 g^{mj}(y) \xi_m \frac{1}{|\xi|_{g(y)}^4}.
    \end{align*}
Due to using geodesic normal coordinates, at the centre $x$ we have $\Gamma_{ij}^k(x)=0$ and $\partial_{x_j}g^{kl}(x)=0$, and hence
\begin{align*}
    r_{-3}(x,\xi)&= 0.
\end{align*}

    For the first term of $T$, namely 
    \[
        T_1:=-2\sum_{i=1}^N \lambda_i C_i M_{e_i} \circ f_i(\Delta_0) \circ M_{e_i},
    \]
    it is clear from the above computations that
    \begin{align*}
        \sigma_{-1}(T_1)(x,\xi) &= 0;\\
        \sigma_{-2}(T_1)(x,\xi)&= -2\sum_{i=1}^N \lambda_i C_i e_i(x)^2\frac{1}{|\xi|_g^2}.
    \end{align*}
    For the second term 
    \[
        T_2:= -(d-2)\sum_{i=1}^N C_i  M_{e_i}\circ f_i(\Delta_0) \circ \nabla_g e_i,
    \]
    we also immediately have
    \begin{align*}
    \sigma_{-1}(T_2)(x,\xi) 
    &= -(d-2) \sum_{i=1}^N C_i \sigma_0(M_{e_i})(x,\xi) \sigma_{-2}(f_i(\Delta_0))(x,\xi) \sigma_1(\nabla_ge_i)(x,\xi) \\
    &= -i(d-2)\sum_{i=1}^N C_i e_i(x) \frac{\xi^\# \cdot \nabla_ge_i}{|\xi|_g^2}.
    \end{align*}
    For $\sigma_{-2}(T_2)$, note similarly to before that
    \begin{align*}
        &\sigma_{-2}(f_i(\Delta_0) \circ \nabla_g e_i)(y,\xi) \\
        &= r_{-3}(y,\xi) ig^{kl}(y)\xi_k \partial_{x_l} e_i(y)-i \sum_{|\alpha|=1} \partial^\alpha_{\xi} \Big(r_{-2}(y,\xi) \Big)\partial^\alpha_{x}\Big(ig^{kl}(y)\xi_k \partial_{x_l} e_i(y)\Big).
    \end{align*}
    Since we are only interested in the value of $\sigma_{-2}(T_2)$ at the centre $x$ where $r_{-3}(x,\xi) =0$, we can ignore the first term. Next, 
    \begin{align*}
        &-i \sum_{|\alpha|=1} \partial^\alpha_{\xi} \Big(r_{-2}(y,\xi) \Big)\partial^\alpha_{x}\Big(ig^{kl}(y)\xi_k \partial_{x_l} e_i(y)\Big)\\
        &= \sum_{j=1}^d \partial_{\xi_j} \Big(\frac{1}{g^{mn}(y)\xi_m\xi_n}\Big)\partial_{x_j}\Big( g^{kl}(y)\xi_k \partial_{x_l} e_i(y)\Big)\\
        &= \sum_{j=1}^d \frac{-2g^{mj}(y)\xi_m}{|\xi|_{g(y)}^4}\bigg(\partial_{x_j}g^{kl}(y) \xi_k \partial_{x_l}e_i(y) + g^{kl}(y) \xi_k \partial_{x_j} \partial_{x_l} e_i(y) \bigg).
    \end{align*}
    At the centre $x$ of the normal coordinates, we therefore have
    \begin{align*}
        \sigma_{-2}(f_i(\Delta_0) \circ \nabla_g e_i)(x,\xi)&= \sum_{j=1}^d \frac{-2g^{mj}(x)\xi_m}{|\xi|_{g}^4}\bigg( g^{kl}(x) \xi_k \partial_{x_j} \partial_{x_l} e_i(x) \bigg)\\
        &= \sum_{j=1}^d \frac{-2\xi_j}{|\xi|_{g}^4}\bigg( \sum_{k=1}^d \xi_k \partial_{x_j} \partial_{x_k} e_i(x) \bigg)\\
        &= -2 \sum_{j,k=1}^d \frac{\xi_j \xi_k \partial_{x_j} \partial_{x_k}e_i(x)}{|\xi|^4_{g}} \\
        &= -2\frac{\nabla^2e_i(\xi,\xi)}{|\xi|_g^4}.
    \end{align*}
    Finally, we obtain
    \begin{align*}
        \sigma_{-2}(T_2)(x,\xi) &= 2(d-2)\sum_{i=1}^N C_i e_i(x)\frac{\nabla^2 e_i(\xi,\xi)}{|\xi|^4_g}.
    \end{align*}
Adding the symbols of $T_1$ and $T_2$ completes the calculation.
\end{proof}

Recall from Lemma~\ref{L:CosphereIntegral} that
\[
\int_{|\xi|_g=1} \nabla^2f(\xi,\xi)\,d\xi = \frac{1}{d} \mathrm{tr}_g \nabla^2 f = -\frac{1}{d}\Delta_g f.
\]
Since we observed that if equation~\eqref{eq:Absurdity2} holds for all $\psi\in C^\infty(M)$, then we must have that $\sigma_{-1}(T)$ and $\sigma_{-2}(T)$ vanish, we would have in particular that
\begin{align*}
\int_{|\xi|_g=1} \sigma_{-2}(T)(x,\xi) \,d\xi &= -2\sum_{i=1}^N \lambda_i C_i e_i(x)^2 +2(d-2)\sum_{i=1}^N C_i e_i(x)\int_{|\xi|_g=1}\nabla^2e_i(\xi,\xi) \, d\xi\\
&= -2\sum_{i=1}^N \lambda_i C_i e_i(x)^2 -2\frac{d-2}{d}\sum_{i=1}^N \lambda_i C_i e_i(x)^2\\
&= -4\frac{d-1}{d} \sum_{i=1}^N \lambda_i C_i e_i(x)^2=0.
\end{align*}
Hence, 
\[
\sum_{i=1}^N \lambda_i C_i e_i(x)^2=0.
\]

We can now iterate all preceding arguments, replacing $C_i$ with $\lambda_i C_i$ in equation~\eqref{eq:Absurdity1}.
\begin{lem}\label{L:WrapUp1}
    Let $(M,g_0)$ be a compact Riemannian manifold  and assume that simultaneously the following to properties hold: 
    \begin{enumerate}
        \item $\lambda_0(g_0) <\cdots <\lambda_N(g_0)$;
        \item for any $\psi \in C^\infty(M)$ giving rise to metrics $g_t = (1+t\psi)g$ and for any $C_i \in \C$ for which
    \[
    \sum_{i=0}^N C_i e_i^2 \equiv 0,
    \]
    there exist analytic functions $C_i(t)$ with $C_i(0)=C_i$ and a small interval $t \in (-\varepsilon, \varepsilon)$ so that
    \[
    \sum_{i=0}^N C_i(t) e_i(t)^2 \equiv 0.
    \]
    \end{enumerate}
    Then the set  of smooth functions  $\{e_i^2 : 0\leq i \leq N\}$ is linearly independent.
\end{lem}
\begin{proof}
    By the previous arguments, if $\sum_{i=0}^N C_i e_i^2\equiv 0$ extends analytically to $ \sum_{i=0}^N C_i(t) e_i(t)^2\equiv 0$ along all perturbations $\psi$ for some small time interval (depending on $\psi$), then
    \[
    \sum_{i=1}^N \lambda_i^k C_i e_i^2 \equiv 0,
    \]
    for all $k \in \N$. 
    Now, since $\lambda_1 < \cdots < \lambda_N$, we have
    \begin{align*}
    0&=\lim_{k\to \infty} \frac{1}{\lambda_N^k} \sum_{i=1}^N \lambda_i ^k C_i e_i^2\\
    &= C_N e_N^2,
    \end{align*}
    hence $C_N =0$. Inductively, $C_i =0 $ for all $1 \leq i \leq N$. The coefficient $C_0$ needs care, since $\lambda_0 = 0$. However, clearly, if $C_0 e_0^2 = 0$, also $C_0 = 0$. This means precisely that $\{e_i^2 : 0\leq i \leq N\}$ is a linearly independent set.
\end{proof}

\begin{proof}[Proof of Theorem~\ref{T:SquareIndependence}]
    Density of the set of metrics for which $ \lambda_0(g) < \cdots < \lambda_N(g)$ and $\{e_i^2:0\leq i \leq N\}$ is linearly independent, follows directly from combining Lemma~\ref{L:Hammer1} and Lemma~\ref{L:WrapUp1}.
    
    We now prove that this set is open. Let $g\in \mathcal{R}(M)$ such that $\lambda_0(g) < \cdots < \lambda_N(g)$ and $\{e_i^2: 0\leq i \leq N\}$ is linearly independent. Since $\lambda_i(g)$ depends continuously on $g\in \mathcal{R}(M)$~\cite{BandoUrakawa1983}, the set of metrics for which $\lambda_0(g) < \cdots < \lambda_N(g)$, is open in $\mathcal{R}(M)$. Then, by applying~\cite[Lemma~3.2]{GreilhuberKepplinger2023}, there exists a neighbourhood $U$ of $g$ in $\mathcal{R}(M)$ and $\delta>0$ such that each $h \in U$ has precisely one eigenvalue in each interval $(\lambda_i-\delta, \lambda_i+\delta)$ with corresponding eigenfunction $e_i(h)$, and there exist continuously differentiable maps 
    \begin{align*}
        F_i: U\subseteq \mathcal{R}(M) &\to C^\infty(M)\\
        h\quad  &\mapsto e_i(h).
    \end{align*}
    Linear independence of $\{e_i^2:0\leq i \leq N\}$ holds if the determinant of the Gram matrix $G = (\langle e_i^2, e_j^2\rangle)_{ij}$ is nonzero. Since this determinant is a polynomial in the coefficients of $G$, it is continuous in $g$. This shows that $\{e_i^2:0\leq i \leq N\}$ stays linearly independent in an open neighbourhood of $g$. 
    \end{proof}
\subsection{Linear independence of auxiliary functions}
Now we will discuss the linear independence of the auxiliary functions
\[
\Big(\frac{2-d}{4}\Delta_0-\lambda_i\Big)e_i^2,
\]
for $1\leq i \leq N$, which appear in the derivative of the eigenvalues in Lemma~\ref{L:EigenvaluePerturbation}. 

\begin{lem}\label{L:Hammer1b}
    Let $M$ be a compact manifold and let $N\in  \N$. 
    Then, the set of metrics,
     which simultaneously satisfy the following three properties, is dense in $\mathcal{R}(M)$:  
    \begin{enumerate}
        \item $\lambda_0(g) < \lambda_1(g)<\cdots<\lambda_N(g)$;
        \item the functions $\{e_i^2: 0\leq i \leq N\}$ are linearly independent;
        \item for any $\psi \in C^\infty(M)$ giving rise to metrics $g_t = (1+t\psi)g$ and for any $C_i \in \C$ for which
        \[
    \sum_{i=1}^N C_i\Big(\frac{2-d}{4}\Delta_0-\lambda_i\Big) e_i^2 \equiv 0,
    \]
    there exist analytic functions $C_i(t)$ with $C_i(0)=C_i$ and a small interval $t \in (-\varepsilon, \varepsilon)$ so that
    \[
    \sum_{i=1}^N C_i(t) \Big(\frac{2-d}{4}\Delta_{g_t}-\lambda_i(t)\Big)e_i(t)^2 \equiv 0.
    \]
    \end{enumerate}
\end{lem}
\begin{proof}
The set of metrics for which $\lambda_0(g) < \lambda_1(g)<\cdots<\lambda_N(g)$ and the functions $\{e_i^2:0\leq i \leq N\}$ are linearly independent, is open and dense due to Theorem~\ref{T:SquareIndependence}. Pick any metric $g_0$ in this set. 
    Along any conformal perturbation $(1+t\psi)g_0$, one can pick eigenfunctions $e_i(t)$ and eigenvalues $\lambda_i(t)$, $0\leq i \leq N$ to depend analytically on $t$ due to Proposition~\ref{P:AnalyticDependence}. Hence also $\big(\frac{2-d}{4}\Delta_{g_t}-\lambda_i(t)\big)e_i(t)^2$ depends analytically on $t$. The proof then continues as that of Lemma~\ref{L:Hammer1}. 
\end{proof}

As before, we want to conclude that each $C_i$ is zero. 
For that purpose, we differentiate (in $t$) the equation
\begin{align*}
    \sum_{i=1}^N C_i(t) \Big(\frac{2-d}{4}\Delta_{g_t}-\lambda_i(t)\Big)e_i(t)^2 \equiv 0.
\end{align*}
Using Lemma~\ref{L:ConformalLaplacian}, we obtain
\begin{align*}
    \sum_{i=1}^N \dot{C}_i\Big(\frac{2-d}{4}\Delta_0-\lambda_i\Big) e_i^2 - C_i\dot{\lambda_i} e_i^2 &= -2\sum_{i=1}^N C_i\Big(\frac{2-d}{4}\Delta_0-\lambda_i\Big) (e_i \dot{e}_i)\\
    &\quad +\frac{2-d}{4} \sum_{i=1}^N C_i \Big(\psi \Delta_0(e_i^2) + (d-2) e_i \nabla_g \psi \cdot \nabla_g e_i \Big).
\end{align*}
Recall from Lemma~\ref{L:EigenfunctionDerivative} that
\[
\dot{e}_i = -\frac{d}{4} \langle e_i^2, \psi\rangle e_i + \lambda_i f_i(\Delta_0)(e_i \psi) + \frac{d-2}{2}f_i(\Delta_0)\big(\nabla_g \psi \cdot \nabla_g e_i\big).
\]
We thus obtain
\begin{align*}
    &\sum_{i=1}^N \Big( \dot{C}_i - \frac{d}{2} C_i\langle e_i^2, \psi\rangle \Big) \Big(\frac{2-d}{4}\Delta_0-\lambda_i\Big) e_i^2 - \sum_{i=1}^N C_i\dot{\lambda_i} e_i^2 \\
    &= -2\sum_{i=1}^N C_i\Big(\frac{2-d}{4}\Delta_0-\lambda_i\Big) \Big( \lambda_i e_i f_i(\Delta_0)(e_i \psi) + \frac{d-2}{2} e_i f_i(\Delta_0)\big(\nabla_g \psi \cdot \nabla_g e_i\big)\Big)\\
    &\quad + \frac{2-d}{4}\sum_{i=1}^N C_i \Big(\psi \Delta_0(e_i^2) + (d-2) e_i \nabla_g \psi \cdot \nabla_g e_i \Big).
\end{align*}

The left-hand side lies in a finite-dimensional space. Hence, for this to be true for all $\psi \in C^\infty(M)$, we need that the principal symbol of the pseudodifferential operator
\begin{align*}
    F &:= -2\sum_{i=1}^N \lambda_i C_i\Big(\frac{2-d}{4}\Delta_0-\lambda_i\Big) \circ   M_{e_i} \circ f_i(\Delta_0) \circ M_{e_i}\\
    &\quad -(d-2)\sum_{i=1}^N C_i \Big(\frac{2-d}{4}\Delta_0-\lambda_i\Big)  \circ M_{e_i} \circ f_i(\Delta_0)\circ \nabla_g e_i\\
    &\quad +\frac{2-d}{4}\sum_{i=1}^N C_i M_{\Delta_0(e_i^2)} -\frac{(d-2)^2}{4} \sum_{i=1}^N C_i M_{e_i} \circ \nabla_{g}e_i
\end{align*}
disappears (here we used Lemma~\ref{L:EigenfunctionDerivative}).

\begin{lem}
Though a priori $F\in \Psi^1_{cl}(M)$, we have that $F \in \Psi^0_{cl}(M)$ with principal symbol
\begin{align*}
    \sigma_0(F)(x,\xi) &= - \frac{(d-2)^2}{2}\sum_{i=1}^N C_i \frac{(\xi \cdot \nabla_g e_i(x))^2}{|\xi|_g^2}+\frac{2-d}{2}\sum_{i=1}^N C_i |\nabla_g e_i|^2.
\end{align*}
\end{lem}
\begin{proof}
Using that 
\[
[\Delta_0, M_{e_i}]f = -2\nabla_g e_i \cdot \nabla_g f + M_{\Delta_0(e_i)}f,
\]
and using that $\Delta_0 \circ f_i(\Delta_0) = 1 + \Psi^{-2}_{cl}(M)$, we can rewrite $F$ modulo order $-1$-operators as
\begin{align*}
    F+ \Psi^{-1}_{cl}(M) &= -\frac{2-d}{2}\sum_{i=1}^N \lambda_i C_i\, \Delta_0\circ   M_{e_i} \circ f_i(\Delta_0) \circ M_{e_i} +\frac{(d-2)^2}{4}\sum_{i=1}^N C_i \, \Delta_0  \circ M_{e_i} \circ f_i(\Delta_0)\circ \nabla_g e_i\\
    &\quad +\frac{2-d}{4}\sum_{i=1}^N C_i M_{\Delta_0(e_i^2)} -\frac{(d-2)^2}{4} \sum_{i=1}^N C_i M_{e_i} \circ \nabla_{g}e_i + \Psi^{-1}_{cl}(M)\\
    &= -\frac{2-d}{2}\sum_{i=1}^N \lambda_i C_i\, M_{e_i^2} +\cancel{\frac{(d-2)^2}{4}\sum_{i=1}^N C_i \,  M_{e_i} \circ \nabla_g e_i} -\frac{(d-2)^2}{2}\sum_{i=1}^N C_i \,  \nabla_g e_i \circ f_i(\Delta_0)\circ \nabla_g e_i\\
    &\quad +\frac{2-d}{4}\sum_{i=1}^N C_i M_{\Delta_0(e_i^2)} -\cancel{\frac{(d-2)^2}{4} \sum_{i=1}^N C_i M_{e_i} \circ \nabla_{g}e_i} + \Psi^{-1}_{cl}(M).
\end{align*}
We indeed see that $F \in \Psi^0_{cl}(M)$, with principal symbol
\begin{align*}
    \sigma_0(F)(x,\xi) &= -\frac{2-d}{2}\sum_{i=1}^N \lambda_i C_i e_i(x)^2 +\frac{(d-2)^2}{2}\sum_{i=1}^N C_i \frac{(\xi \cdot \nabla_g e_i(x))^2}{|\xi|_g^2} +\frac{2-d}{4}\sum_{i=1}^N C_i \Delta_0(e_i^2)\\
    &= \frac{(d-2)^2}{2}\sum_{i=1}^N C_i \frac{(\xi \cdot \nabla_g e_i(x))^2}{|\xi|_g^2}+\frac{d-2}{2}\sum_{i=1}^N C_i |\nabla_g e_i|^2.
\end{align*}
\end{proof}
Recall from Lemma~\ref{L:CosphereIntegral} that
\[
\int_{|\xi|_g=1} (\xi \cdot \nabla_g f)^2 \, d\xi = \frac{1}{d} |\nabla_g f|^2.
\]
Therefore, if the principal symbol of $F$ vanishes, then by taking the fibre-wise integral over the cosphere bundle $S^*_xM$ we obtain that
\begin{align*}
0&=\int_{|\xi|_g=1} \sigma_0(F)(x,\xi) \, d\xi\\
&=  \frac{(d-2)^2}{2d}\sum_{i=1}^N C_i |\nabla_g e_i|^2+\frac{d-2}{2}\sum_{i=1}^N C_i |\nabla_g e_i|^2\\
&= \frac{(d-1)(d-2)}{d}\sum_{i=1}^N C_i |\nabla_g e_i|^2.
\end{align*}
For $d\geq 3$ this of course means that $\sum_{i=1}^N C_i |\nabla_g e_i|^2$. But, that finally implies that
\begin{align*}
    0&\equiv \sum_{i=1}^N C_i \Big(\frac{2-d}{4}\Delta_0-\lambda_i\Big)e_i^2\\
    &= \frac{d-2}{2} \sum_{i=1}^N C_i |\nabla_g e_i|^2 - \frac{d}{2}\sum_{i=1}^N \lambda_i C_i e_i^2 \\
    &= - \frac{d}{2}\sum_{i=1}^N \lambda_i C_i e_i^2.
\end{align*}
By the linear independence of $e_i^2$, $1\leq i \leq N$ due to  Theorem~\ref{T:SquareIndependence}, we obtain that $C_i =0$ for all $1 \leq i \leq N$. In summary, we have shown the following.
\begin{lem}\label{L:WrapUp1b}
    If $(M,g_0)$ is a compact Riemannian manifold of dimension $d\geq 3$ such that 
    \begin{enumerate}
        \item  $\lambda_0(g_0) <\cdots <\lambda_N(g_0)$;
        \item  $\{ e_i(g_0)^2: 0\leq i \leq N\}$ is a linearly independent set;
        \item along any conformal perturbation $g_t = (1+t\psi) g_0$ and for all $C_i\in \C$ for which
    \[
    \sum_{i=1}^N C_i \Big(\frac{2-d}{4}\Delta_0-\lambda_i\Big) e_i^2\equiv 0,
    \]
    there exist real-analytic functions $C_i(t)$ with $C_i(0)=C_i$ for which on a small time interval
    \[
    \sum_{i=1}^N C_i(t) \Big(\frac{2-d}{4}\Delta_{g_t}-\lambda_i(t)\Big) e_i(t)^2\equiv 0,
    \]
    \end{enumerate}
    then we must have that the functions 
\[
\{\Big(\frac{2-d}{4}\Delta_0-\lambda_i\Big) e_i^2: 1 \leq i \leq N\},
\]
are linearly independent.
\end{lem}

We finally conclude with the main result of this subsection. Note that this theorem is different from Theorem~\ref{T:SquareIndependence} in that we do not show openness. It is unclear to the authors whether the dependency of $(\frac{2-d}{4}\Delta_g - \lambda_i(g)) e_i^2(g)$ on $g$ is continuous. We furthermore do not need this for our other results.

\begin{thm}\label{T:CursedIndependence}
    Let $M$ be a smooth compact manifold. For each $N\in \N$, the set of metrics $g$ in $\mathcal{R}(M)$ for which simultaneously
    \begin{enumerate}
        \item $\lambda_0(g) < \lambda_1(g)<\cdots<\lambda_N(g) $;
        \item  $\{e_i^2 : 0\leq i \leq N\}$ is a set of linearly independent functions;
        \item  $\{(\frac{2-d}{4}\Delta_0 - \lambda_i) e_i^2 : 1 \leq i \leq N \}$ is a set of linearly independent functions,
    \end{enumerate}
     is dense in $\mathcal{R}(M)$.
\end{thm}
\begin{proof}
    For $d=2$, we have $(\frac{2-d}{4}\Delta_0 - \lambda_i) e_i^2 = -\lambda_i e_i^2$, in which case this theorem reduces to Theorem~\ref{T:SquareIndependence}. For $d\geq 3$, this follows from combining Lemma~\ref{L:Hammer1b} and Lemma~\ref{L:WrapUp1b}.
\end{proof}

\subsection{Linear independence of spectrum over the rationals}
We first of all note that the spectrum of $\Delta_g$ (omitting $\lambda_0=0$) is linearly independent over $\mathbb{Q}$ if and only if it is linearly independent over $\mathbb{Z}$. Hence, to prove that the set of metrics for which the spectrum of $\Delta_g$ is linearly independent over $\mathbb{Q}$ is residual, it suffices to prove that the set of metrics for which $\lambda_1(g), \ldots, \lambda_N(g)$ are linearly independent over $\mathbb{Z}\cap[-K,K]$ is open and dense for all $N,K \in \N$.

This section is an adaptation of the strategy in~\cite{BleeckerWilson1980,BandoUrakawa1983}. To prove that the set of metrics with no linear dependencies in the spectrum $\leq N$ is dense, we only need Theorem~\ref{T:SquareIndependence} for manifolds of dimension $d=2$. For $d\geq 3$, we need the full power of Theorem~\ref{T:CursedIndependence}.

\begin{thm}\label{T:NoResonances}
    Let $M$ be a smooth compact manifold. For each $N,K \in \N$, the set of metrics $g$ in $\mathcal{R}(M)$ for which $\lambda_1(g),\ldots, \lambda_N(g)$ are linearly independent over $\mathbb{Z} \cap [-K,K]$, is open and dense.
\end{thm}
\begin{proof}
    The condition that $\lambda_1(g),\ldots, \lambda_N(g)$ are linearly independent over $\mathbb{Z} \cap [-K,K]$ can be expressed as
    \[
    F_{k_1,\ldots,k_N}(g) := \sum_{i=1}^N k_i \lambda_i(g) \neq 0,
    \]
    for all such defined functions, with $k_i \in \mathbb{Z}\cap [-K,K]$. This is, of course, a finite family of functions. Since the eigenvalues continuously on $g\in \mathcal{R}(M)$~\cite[Theorem~2.2]{BandoUrakawa1983}, the functions $F_{k_1,\ldots,k_N}$ are continuous. Hence it is clear that the set of metrics for which each $F_{k_1,\ldots,k_N}$ is non-zero is open.
    
    Now we prove density. Take any metric $g_0 \in \mathcal{R}(M)$ such that both
    \begin{enumerate}
        \item $\lambda_0(g) < \lambda_1(g)<\cdots<\lambda_N(g) $;
        \item  $\{(\frac{2-d}{4}\Delta_0 - \lambda_i) e_i^2 : 1 \leq i \leq N \}$ is a set of linearly independent functions.
    \end{enumerate}
    These metrics are dense in $\mathcal{R}(M)$ according to Theorem~\ref{T:CursedIndependence}. 
    Perturbing this metric to $g_t = (1+t\psi)g_0$ with $\psi \in C^\infty(M)$, we have for $\lambda_i(t)$, $1 \leq i \leq N$, that
    \[
    \lambda_i(t) = \lambda_i + t \dot{\lambda_i} + t^2 \mu_i(t),
    \]
    where $\mu_i(t)$ is an analytic function. Furthermore, by Lemma~\ref{L:EigenvaluePerturbation}, we obtained
\[
    \dot{\lambda_i} =  \big\langle \psi,  \Big(\frac{2-d}{4}\Delta_0-\lambda_i\Big)e_i^2\big\rangle_0.
\]
Hence, we can choose $\psi \in C^\infty(M)$ to set the values of $\dot{\lambda_i}$ for $1\leq i \leq N$ arbitrarily by Lemma~\ref{L:LinIndep}. If $F_{k_1,\ldots,k_N}(g) = 0$ for some choices of $k_1,\ldots, k_N$, we can find $\psi$ along which $\sum_{i=1}^N k_i \dot{\lambda_i}\neq 0$ for each such relation. Therefore, for some interval $0\neq t\in (-\varepsilon, \varepsilon)$, we obtain that $F_{k_1,\ldots,k_N}(g_t)\neq 0$. This implies the claim of the theorem.
\end{proof}

\begin{rem}\label{R:Gradients}
    The proof of Theorem~\ref{T:NoResonances} shows that, for a dense set of metrics $g_0\in \mathcal{R}(M)$ for which $\lambda_0<\ldots<\lambda_N$, we can find for any $(C_1,\ldots,C_N)\in \R^N$ a conformal perturbation $\psi \in C^\infty(M)$ so that, along $g_t=(1+t\psi)g_0$, we have
    \[
    \dot{\lambda_i} = C_i, \quad i=1,\ldots,N.
    \]
    This illustrates that for any finite family of $C^1$-functions $\{F_k:\R^{N}\to \C\}_{k=1}^K$ with $\nabla F_k$ nowhere zero, the set of metrics $g \in \mathcal{R}(M)$ such $\lambda_1(g)<\cdots < \lambda_N(g)$ and for all $k=1,\ldots, K$
    \[
    F_k(\lambda_1(g),\ldots, \lambda_{N_k}(g)) \neq 0,
    \]
    is open and dense. The same holds for any finite family of non-zero real-analytic functions $\{F_k:\R^{N}\to \C\}_{k=1}^K$. Namely, at a zero of order $n$, there exists $v \in \mathbb{R}^N$ with $D^n F_k(v,\ldots, v) \neq 0$. Choosing $\psi \in C^\infty(M)$ such that $\dot{\lambda_i}=v_i$, $1 \leq i \leq N$, provides that 
    \[
    F_k(\lambda_1(t),\ldots, \lambda_N(t))= \frac{t^k}{k!} D^n F_k(v,\ldots, v) + O(t^{k+1}),
    \]
    which is nonzero in some punctured interval $(-\varepsilon, \varepsilon)\setminus\{0\}$.
\end{rem}

\subsection{Linear independence of products}
Recall the multi-index notation
we write
\[
  e^\alpha := \prod_{i\ge 1} e_i^{\alpha_i},\qquad
  |\alpha| := \sum_{i\ge1}\alpha_i,
\]
for finitely supported multi-indices $\alpha=(\alpha_i)_{i\ge1}\in\mathbb N^{(\infty)}$,
with the convention $e^0 = e_0$. 
We now perform the same argument as in Section~\ref{SS:Squares}, but with the relations
\[
\sum_{\substack{\alpha \in\mathbb {N}^{N}\\ |\alpha|\leq L}}C_{\alpha} e^\alpha \equiv 0.
\]

\begin{lem}\label{L:Hammer2}
    Let $M$ be a smooth compact manifold and let $N,K, L\in  \N$. Then, the set of metrics for which simultaneously
    \begin{enumerate}
        \item $\lambda_1(g), \ldots, \lambda_N(g)$ are linearly independent over $\mathbb{Z}\cap[-K,K]$;
        \item for all $\psi \in C^\infty(M)$ giving rise to metrics $g_t$ and for all $C_{i_1\ldots i_n} \in \C$ for which
    \[
    \sum_{\substack{\alpha \in\mathbb {N}^{N}\\ |\alpha|\leq L}}C_{\alpha} e^\alpha \equiv 0,
    \]
    there exist analytic functions $C_{\alpha}(t)$ with $C_\alpha(0)=C_\alpha$ and for which along some small interval $t \in (-\varepsilon, \varepsilon)$
    \[
    \sum_{\substack{\alpha \in\mathbb {N}^{N}\\ |\alpha|\leq L}}C_{\alpha}(t) e^\alpha(t) \equiv 0,
    \]
    \end{enumerate} 
    is dense in $\mathcal{R}(M)$. 
\end{lem}
\begin{proof}
    Follows from the same proof as Lemma~\ref{L:Hammer1}, using the Gram matrix for the functions $\{e^\alpha : \alpha \in\mathbb {N}^{N}, \, |\alpha|\leq L\}$, and now using that the set of metrics for which $\lambda_1(g), \ldots, \lambda_N(g)$ are linearly independent over $\mathbb{Z}\cap[-K,K]$ is dense in $\mathcal{R}(M)$ due to Theorem~\ref{T:NoResonances}.
\end{proof}

Now, if we differentiate
\[
     \sum_{\substack{\alpha \in\mathbb {N}^{N}\\ |\alpha|\leq L}}C_{\alpha}(t) e^\alpha(t) \equiv 0,
    \]
we obtain
\begin{align*}
 \sum_{\substack{\alpha \in\mathbb {N}^{N}\\ |\alpha|\leq L}}\dot{C}_{\alpha} e^\alpha &= - \sum_{\substack{\alpha \in\mathbb {N}^{N}\\ |\alpha|\leq L}}C_{\alpha}  \Big(\sum_{j=1}^N e_{1}^{\alpha_1} \cdots (\alpha_j e_j^{\alpha_j-1} \dot{e}_{j})\cdots e_{N}^{\alpha_N}\Big).
\end{align*}
By Lemma~\ref{L:EigenfunctionDerivative}, this means that
\begin{equation}\label{eq:Absurdity3}\begin{split}
&\sum_{\substack{\alpha \in\mathbb {N}^{N}\\ |\alpha|\leq L}}\Big(\dot{C}_{\alpha} - \frac{d}{4}C_{\alpha} \langle \Big(\sum_{j=1}^N \alpha_j e_{j}^2\Big),\psi\rangle_0 \Big)  e^\alpha \\
&= -\sum_{\substack{\alpha \in\mathbb {N}^{N}\\ |\alpha|\leq L}}C_{\alpha} \Big(\sum_{j=1}^N \alpha_j e_{1}^{\alpha_1} \cdots e_{j}^{\alpha_j-1}\big(\lambda_{j} f_{j}(\Delta_0)(e_{j}\psi)\big)\cdots e_{N}^{\alpha_N}\Big)\\
&\quad - \frac{d-2}{2}\sum_{\substack{\alpha \in\mathbb {N}^{N}\\ |\alpha|\leq L}}C_{\alpha} \Big(\sum_{j=1}^N \alpha_j e_{1}^{\alpha_1} \cdots e_{j}^{\alpha_j-1}\big( f_{j}(\Delta_0)(\nabla_g\psi \cdot \nabla_g e_{j})\big)\cdots e_{N}^{\alpha_N}\Big)
\end{split}
\end{equation}
As argued in Section~\ref{SS:Squares}, the left-hand side lies in a finite-dimensional space spanned by $
\{e^\alpha: \alpha \in\mathbb {N}^{N},\ |\alpha|\leq L\}$, and the right-hand side is a pseudodifferential operator applied to $\psi$, namely
\begin{equation}\label{eq:PsDOS}\begin{split}
    S&:=-\sum_{\substack{\alpha \in\mathbb {N}^{N}\\ |\alpha|\leq L}}\sum_{j=1}^N \alpha_j C_{\alpha} M_{e_{1}^{\alpha_1} \cdots e_{j}^{\alpha_j-1}\cdots e_{N}^{\alpha_N}} \circ \lambda_{j} f_{j}(\Delta_0) \circ M_{e_{j}}\\
&\quad - \frac{d-2}{2} \sum_{\substack{\alpha \in\mathbb {N}^{N}\\ |\alpha|\leq L}}\sum_{j=1}^N \alpha_j C_{\alpha} M_{e_{1}^{\alpha_1} \cdots e_{j}^{\alpha_j-1}\cdots e_{N}^{\alpha_N}} \circ  f_{j}(\Delta_0) \circ \nabla_g e_{j}.
\end{split}
\end{equation}

From the computations in Lemma~\ref{L:TPrincipalSymb}, we obtain the following.
\begin{lem}\label{L:SPrincipalSymb}
    For $S \in \Psi^{-1}_{cl}(M)$ defined in equation~\eqref{eq:PsDOS}, we have for $x$ in the centre of geodesic normal coordinates,
        \begin{align*}
        \sigma_{-1}(S)(x,\xi) &= -i\frac{d-2}{2} \sum_{\substack{\alpha \in\mathbb {N}^{N}\\ |\alpha|\leq L}}\sum_{j=1}^N \alpha_j C_{\alpha} e_{1}(x)^{\alpha_1} \cdots e_{j}(x)^{\alpha_j-1} \frac{\xi \cdot \nabla_g e_{j}}{|\xi|_g^2} \cdots e_{N}(x)^{\alpha_N};\\
        \sigma_{-2}(S)(x,\xi)&= -\sum_{\substack{\alpha \in\mathbb {N}^{N}\\ |\alpha|\leq L}}  C_{\alpha} \Big(\sum_{j=1}^N \alpha_j\lambda_{j}\Big)e^\alpha(x) \frac{1}{|\xi|_g^2}\\
        &\quad +(d-2) \sum_{\substack{\alpha \in\mathbb {N}^{N}\\ |\alpha|\leq L}}\sum_{j=1}^N \alpha_j C_{\alpha} e_{1}(x)^{\alpha_1} \cdots e_{j}(x)^{\alpha_j-1} \frac{\nabla^2 e_{j}(\xi,\xi)}{|\xi|_g^4} \cdots e_{N}(x)^{\alpha_N}.
    \end{align*}
\end{lem}

For equation~\eqref{eq:Absurdity3} to hold for all $\psi \in C^\infty(M)$, we must have that $\sigma_{-1}(S)$ and $\sigma_{-2}(S)$ vanish. 
As before, this implies that the fibre-wise integral of $\sigma_{-2}(S)$ over the cosphere bundle $S_x^*M$ vanishes,
\begin{align*}
0&=\int_{|\xi|_g=1} \sigma_{-2}(S)(x,\xi) \,d\xi \\
&= -\sum_{\substack{\alpha \in\mathbb {N}^{N}\\ |\alpha|\leq L}}  C_{\alpha} \Big(\sum_{j=1}^N \alpha_j\lambda_{j}\Big)e^\alpha(x)\\
&\quad +(d-2) \sum_{\substack{\alpha \in\mathbb {N}^{N}\\ |\alpha|\leq L}}\sum_{j=1}^N \alpha_j C_{\alpha} e_{1}(x)^{\alpha_1} \cdots e_{j}(x)^{\alpha_j-1} \Big(\int_{|\xi|_g=1}\nabla^2 e_{j}(\xi,\xi)\, d\xi\Big) \cdots e_{N}(x)^{\alpha_N}\\
&= -2 \frac{d-1}{d} \sum_{\substack{\alpha \in\mathbb {N}^{N}\\ |\alpha|\leq L}}  C_{\alpha} \Big(\sum_{j=1}^N \alpha_j\lambda_{j}\Big)e^\alpha(x).
\end{align*}

\begin{lem}\label{L:WrapUp2}
    Let $N,L,K \in \N$ with $K \geq L$. If $(M,g_0)$ is a compact Riemannian manifold such that 
    \begin{enumerate}
        \item  $\lambda_1, \ldots, \lambda_N$ are linearly independent over $\mathbb{Z} \cap [-K,K]$;
        \item for all $\psi \in C^\infty(M)$ giving rise to metrics $g_t$ and for all $C_{\alpha} \in \C$ for which
    \[
    \sum_{\substack{\alpha \in\mathbb {N}^{N}\\ |\alpha|\leq L}}C_{\alpha} e^\alpha \equiv 0,
    \]
    there exist real-analytic functions $C_{\alpha}(t)$ with $C_\alpha(0)=C_\alpha$ for which along some small interval $t \in (-\varepsilon, \varepsilon)$
    \[
    \sum_{\substack{\alpha \in\mathbb {N}^{N}\\ |\alpha|\leq L}}C_{\alpha}(t) e^\alpha(t) \equiv 0,
    \]
    \end{enumerate}
    then we must have that $(e^\alpha)_{\alpha \in\mathbb {N}^{N},\ |\alpha|\leq L}$ is a linearly independent family.
\end{lem}
\begin{proof}
    The property that $\lambda_1, \ldots, \lambda_N$ are linearly independent over $\mathbb{Z} \cap [-K,K]$, with $K \geq L$, implies that for $\alpha, \beta \in\mathbb {N}^{N}$ with $|\alpha|, |\beta| \leq L$,
    \[
    \sum_{j=1}^N \alpha_j\lambda_{j} = \sum_{j=1}^N \beta_j\lambda_{j} \quad \iff \quad \alpha = \beta.
    \]
    Hence, each such sum takes a unique value, and the result then follows from the same proof as Lemma~\ref{L:WrapUp1}.
\end{proof}

We finally obtain our main theorem.
\begin{thm}\label{T:FinalMain}
    Let $M$ be a smooth compact manifold. For each $N,K,L \in \N$ with $K \geq L$, the set of metrics $g \in \mathcal{R}(M)$ for which simultaneously
    \begin{enumerate}
        \item $\lambda_1(g), \ldots, \lambda_N(g)$ are linearly independent over $\mathbb{Z}\cap[-K,K]$;
        \item the products of eigenfunctions $\{e^\alpha: \alpha \in\mathbb {N}^{N},\ |\alpha|\leq L\}$ are linearly independent,
    \end{enumerate} is open and dense.
\end{thm}
\begin{proof}
    That this set of metrics is open, follows from the same argument as in the proof of Theorem~\ref{T:SquareIndependence}, after noting that by~\cite[Lemma~3.2]{GreilhuberKepplinger2023}, in a neighbourhood of $g$ for which $\lambda_1(g),\ldots,\lambda_N(g)$ are linearly independent over $\mathbb{Z}\cap[-K,K]$, there exists a continuous function
    \begin{align*}
        F: U\subseteq \mathcal{R}(M) &\to C^\infty(M)\\
        h \quad & \mapsto e_i(h).
    \end{align*}
    
    Density follows from combining Lemma~\ref{L:Hammer2} and Lemma~\ref{L:WrapUp2}.
\end{proof}

\begin{rem}\label{R:FinalGradients}
    Following Remark~\ref{R:Gradients}, we can generalise the linear independence of $\lambda_1(g),\ldots, \lambda_N(g)$ over $\mathbb{Z}\cap [-K,K]$ in Theorem~\ref{T:FinalMain} to the property that 
    \[
    F_k(\lambda_1(g),\ldots,\lambda_N(g)) \neq 0, \quad 1\leq k \leq M,
    \]
    for some finite collection of $C^1$-functions $\{F_k:\mathbb{R}^{N} \to \mathbb{C}\}_{k=1}^M$ with $\nabla F_k$ nowhere zero. Or likewise, for a finite collection of non-zero real-analytic functions $\{F_k:\mathbb{R}^{N} \to \mathbb{C}\}_{k=1}^M$. This provides the claim Remark~\ref{R:Strongest}. Taking the functions
    \[
    F_k(\lambda_1(g),\ldots, \lambda_{N_k}(g)) = \sum_{j=1}^{N_k}\sum_{l=1}^{M_k} m_{j,k_l} \lambda_j(g)^{q_{j,k_l}}, \quad m_{j,k_l}\in \mathbb{Z},q_{j,k_l} \in \mathbb{Q},
    \]
    in some suitable enumeration provides the statement as given in the~\autoref{T:MainThm}.
\end{rem}

\bibliographystyle{abbrv}
\bibliography{bibliography}

\end{document}